\documentclass{amsart}

\usepackage{xypic}
\usepackage{amsrefs}
\usepackage{longtable,booktabs}
\usepackage{xcolor}
\usepackage[pdf]{pstricks}
\usepackage{pst-plot}
\usepackage{lscape}

\newtheorem{thm}{Theorem}[section]
\newtheorem{conj}[thm]{Conjecture}

\theoremstyle{definition} 
\newtheorem{remark}[thm]{Remark}

\newcommand{\C}{\mathbb{C}}	
\newcommand{\F}{\mathbb{F}} 
\newcommand{\Z}{\mathbb{Z}} 

\newcommand{\map}{\rightarrow}

\DeclareMathOperator{\Ext}{Ext}
\DeclareMathOperator{\Sq}{Sq}
\DeclareMathOperator{\Hom}{Hom}

\newcommand{\mydot}{\mathord{\cdot}}

\newgray{gridline}{0.8}

\newrgbcolor{vzerocolor}{1 0 0}
\newrgbcolor{vonecolor}{0.1 0.1 0.7}
\newrgbcolor{vtwocolor}{0.1 0.8 0.1}
\newrgbcolor{vthreecolor}{0.5 0.5 0.5}

\newcommand{\cirrad}{0.1}

\begin{document}

\title{Stable homotopy groups of spheres}

\author{Daniel C.\ Isaksen}
\address{Department of Mathematics\\
Wayne State University\\
Detroit, MI 48202, USA}
\email{isaksen@wayne.edu}

\thanks{The first author was supported by NSF grant
DMS-1606290.
The second author was supported by grant NSFC-11801082. 
The third author was supported by NSF grant DMS-1810638.
Many of the associated machine computations were performed on the
Wayne State University Grid high performance computing cluster.}

\author{Guozhen Wang}
\address{Shanghai Center for Mathematical Sciences, Fudan University, Shanghai, China, 200433}
\email{wangguozhen@fudan.edu.cn}

\author{Zhouli Xu}
\address{Department of Mathematics, Massachusetts Institute of Technology, Cambridge, MA 02139}
\email{xuzhouli@mit.edu}



\date{\today}

\begin{abstract}
We discuss the current state of knowledge of stable homotopy groups
of spheres.  We describe a new computational method
that yields a streamlined computation of the first
61 stable homotopy groups, and gives new information
about the stable homotopy groups in dimensions 62 through 90.
The method relies more heavily on machine computations than
previous methods, and is therefore less prone to error.
The main mathematical tool is the Adams spectral sequence.
\end{abstract}

\maketitle

\section{Background and history}
\label{sctn:background}

One of the most fundamental problems in topology is to determine
the set of homotopy classes $[S^{n+k}, S^n]$
of continuous based maps $f: S^{n+k} \map S^n$ between spheres.
For $n+k \geq 1$, these sets have a natural group structure,
and they are abelian when $n+k \geq 2$.  Despite their
essential topological importance, they are notoriously 
difficult to compute.  Detailed computations exist only
when $k$ is at most approximately $30$
\cite{Mimura64} \cite{MMO75} \cite{MT63} \cite{Oda77}
\cite{Ravenel86}*{Table A3.6} \cite{Toda62}.

Theorem \ref{thm:intro} summarizes some basic structural results
about $[S^{n+k}, S^n]$.
\begin{thm}
\label{thm:intro}
\mbox{}
\begin{enumerate}
\item
$[S^{n+k}, S^n] = 0$ when $k < 0$.
\item
\label{part:0-stem}
$[S^n, S^n] = \Z$.
\item
\label{part:finite}
$[S^{n+k}, S^n]$ is a finite group, except
when $k = 0$, or $n$ is even and $k = n-1$
(Serre finiteness theorem \cite{Serre53}).
\end{enumerate}
\end{thm}

Unfortunately, Theorem \ref{thm:intro} gives little information
about the groups $[S^{n+k}, S^n]$ when $k > 0$.

The Freudenthal Suspension Theorem \cite{Freudenthal37} provides
a relationship between the groups $[S^{n+k}, S^n]$
for fixed $k$ and varying $n$.  The suspension map induces a
sequence
\[
\cdots \map [S^{n-1+k}, S^{n-1}] \map [S^{n+k}, S^n] \map
[S^{n+1+k}, S^{n+1}] \map \cdots
\]
of group homomorphisms, and these homomorphisms
are in fact isomorphisms when $n > k + 1$.
The ``stable value" $[S^{n+k}, S^n]$ for $n$ sufficiently large
is known as the $k$th stable homotopy group $\pi_k$.

The stable homotopy groups $\pi_k$ enjoy
additional structure that make them more amenable to computation
than the unstable groups $[S^{n+k}, S^n]$.
The rest of this article is entirely concerned with the stable
homotopy groups.  While the study of stable homotopy groups gives
much information about the structure of $[S^{n+k}, S^n]$,
it does not give complete information.  There are
still plenty of groups for which $n \leq k+1$, and stable homotopy
group information does not tell us much about those groups.

Theorem \ref{thm:intro}
implies that
$\pi_k$ is zero if $k < 0$;
that $\pi_0$ is isomorphic to $\Z$; and that
$\pi_k$ is a finite group for all $k > 0$.
Because each group $\pi_k$ is finite, it makes sense to study the
groups one prime at a time.  More specifically, we can
compute the $p$-primary component of the group $\pi_k$ 
for all primes $p$,
and then reassemble these components into a uniquely
determined finite abelian group.

History has demonstrated the effectiveness of this $p$-primary
approach.  For the remainder of this article, we will focus
on the $2$-primary components of $\pi_k$,
except that
we record odd primary information in Table \ref{tab:main}
as a convenience for the reader.
We write $\pi_k^{\wedge}$
for the $2$-primary component of $\pi_k$.  There are plenty of interesting
phenomena to study at odd primes, but we leave that discussion
for other work.  
See Ravenel's comprehensive book \cite{Ravenel86} 
for an extensive source of information on
computations at odd primes.

We will not discuss the long history of stable homotopy group
computations thoroughly. 
  Some of the key results include
the original definitions of homotopy groups by {\v C}ech and Hurewicz;
the discovery of the Hopf maps \cite{Hopf31};
the connection to framed cobordism due to Pontryagin \cite{Pontryagin50}
\cite{Rohlin51};
Serre's method using fiber sequences and
Eilenberg-Mac Lane spaces \cite{Serre51};
the Adams spectral sequence \cite{Adams58};
the May spectral sequence \cite{May64};
Toda's work involving higher compositions and the EHP sequence
\cite{Toda62}; deeper analysis of the Adams spectral sequence
by Mahowald (with Barratt, Bruner, and Tangora) \cite{BMT70} \cite{Bruner84} \cite{MT67};
the Adams-Novikov spectral sequence \cite{Novikov67};
and Kochman's approach involving
the Atiyah-Hirzebruch spectral sequence for the Brown-Peterson
spectrum $BP$ \cite{Kochman90} \cite{KM93}.

The stable homotopy groups have important applications in the study
of high-dimensional manifolds.  
For general values of $n$, Milnor showed that 
the topological manifold $S^n$
may possess more than one non-diffeomorphic smooth structure
\cite{Milnor56}.  
Kervaire and Milnor reduced 
the classification of smooth structures to
a computation of stable homotopy groups
\cite{KM77}.  More precisely, this classification
relies on the $v_1$-torsion subgroups of $\pi_k$
defined at the end of Section \ref{sctn:v1},
as well as the
Kervaire invariant.  Despite the recent breakthrough 
of Hill, Hopkins, and Ravenel \cite{HHR16},
there is still one remaining unknown value of the Kervaire
invariant in dimension $126$.
Previous work of the authors established unique
smooth structures on the spheres $S^{56}$ and $S^{61}$.
See the introduction of \cite{WX17} for a more detailed
discussion of these ideas.

\section{Stable homotopy group computations}

We use the $\C$-motivic homotopy theory of Morel and Voevodsky \cite{MV99},
which has a richer
structure than classical homotopy theory, to deduce new
information about stable homotopy groups.
In practice, our procedure works remarkably well.  Already
we have obtained nearly complete information to dimension 90,
extending well beyond the previously known range that
ended at dimension 61.
In the history of stable homotopy group computations, 
there are only a few occurrences of breakthroughs of a similar
magnitude, including work of Hopf, Serre, Adams, May, Mahowald,
and others.

\begin{thm}
\label{thm:main}
Table \ref{tab:main} describes the stable homotopy groups $\pi_k$
for all values of $k$ up to $90$.
\end{thm}

We adopt the following
notation in Table \ref{tab:main}.  
An integer $n$ stands for the cyclic abelian group $\Z/n$;
the expression $n\mydot m$ stands for the direct sum $\Z/n \oplus \Z/m$;
and $n^j$ stands for the direct sum of $j$ copies of $\Z/n$.
The horizontal line after dimension 61 indicates the range
in which our computations are new information.

Table \ref{tab:main} describes each group $\pi_k$ as the direct sum
of three subgroups.
The first and third columns describe the 2-primary component of each group,
while the second and third columns describe the odd primary components of each group.
See Section \ref{sctn:v1} below for a brief explanation of the meaning
of the $v_1$-torsion and $v_1$-periodic subgroups.

There remain some uncertainties in the $v_1$-torsion subgroups at the
prime $2$.  The first such uncertainties occur in dimensions $69$ and $70$,
and there are additional uncertainties beyond dimension $80$.
In most cases, these uncertainties mean that the order of $\pi_k$ is 
known only up to a factor of $2$.  In a few cases, the additive group
structure is also undetermined.
See Section \ref{sctn:unknown} below for more discussion.

\begin{longtable}{l|lll}
\caption{Stable homotopy groups up to dimension $90$
\label{tab:main} 
} \\
\toprule
$k$ & $v_1$-torsion & $v_1$-torsion & $v_1$-periodic \\
& at the prime 2 & at odd primes & \\
\midrule \endfirsthead
\caption[]{Stable homotopy groups up to dimension $90$} \\
\toprule
$k$ & $v_1$-torsion & $v_1$-torsion & $v_1$-periodic \\
& at the prime 2 & at odd primes & \\
\midrule \endhead
\bottomrule \endfoot
$1$ & $0$ & $0$ & $2$ \\
$2$ & $0$ & $0$ & $2$ \\
$3$ & $0$ & $0$ & $8\mydot3$ \\
$4$ & $0$ & $0$ & $0$ \\
$5$ & $0$ & $0$ & $0$ \\
$6$ & $2$ & $0$ & $0$ \\
$7$ & $0$ & $0$ & $16\mydot3\mydot5$ \\
$8$ & $2$ & $0$ & $2$\\
$9$ & $2$ & $0$ & $2^2$ \\
$10$ & $0$ & $3$ & $2$ \\
$11$ & $0$ & $0$ & $8\mydot9\mydot7$ \\
$12$ & $0$ & $0$ & $0$ \\
$13$ & $0$ & $3$ & $0$ \\
$14$ & $2^2$ & $0$ & $0$ \\
$15$ & $2$ & $0$ & $32\mydot3\mydot5$ \\
$16$ & $2$ & $0$ & $2$ \\
$17$ & $2^2$ & $0$ & $2^2$ \\
$18$ & $8$ & $0$ & $2$ \\
$19$ & $2$ & $0$ & $8\mydot3\mydot11$ \\
$20$ & $8$ & $3$ & $0$ \\
$21$ & $2^2$ & $0$ & $0$ \\
$22$ & $2^2$ & $0$ & $0$ \\
$23$ & $2\mydot8$ & $3$ & $16\mydot9\mydot5\mydot7\mydot13$ \\
$24$ & $2$ & $0$ & $2$ \\
$25$ & $0$ & $0$ & $2^2$ \\
$26$ & $2$ & $3$ & $2$ \\
$27$ & $0$ & $0$ & $8\mydot3$ \\
$28$ & $2$ & $0$ & $0$ \\
$29$ & $0$ & $3$ & $0$ \\
$30$ & $2$ & $3$ & $0$ \\
$31$ & $2^2$ & $0$ & $64\mydot3\mydot5\mydot17$ \\
$32$ & $2^3$ & $0$ & $2$ \\
$33$ & $2^3$ & $0$ & $2^2$ \\
$34$ & $4\mydot2^2$ & $0$ & $2$ \\
$35$ & $2^2$ & $0$ & $8\mydot27\mydot7\mydot19$ \\
$36$ & $2$ & $3$ & $0$ \\
$37$ & $2^2$ & $3$ & $0$ \\
$38$ & $2\mydot4$ & $3\mydot5$ & $0$ \\
$39$ & $2^5$ & $3$ & $16\mydot3\mydot25\mydot11$ \\
$40$ & $2^3\mydot4$ & $3$ & $2$ \\
$41$ & $2^3$ & $0$ & $2^2$ \\
$42$ & $2\mydot8$ & $3$ & $2$ \\
$43$ & $0$ & $0$ & $8\mydot3\mydot23$ \\
$44$ & $8$ & $0$ & $0$ \\
$45$ & $2^3\mydot16$ & $9\mydot5$ & $0$ \\
$46$ & $2^4$ & $3$ & $0$ \\
$47$ & $2^3\mydot4$ & $3$ & $32\mydot9\mydot5\mydot7\mydot13$ \\
$48$ & $2^3\mydot4$ & $0$ & $2$ \\
$49$ & $0$ & $3$ & $2^2$ \\
$50$ & $2^2$ & $3$ & $2$ \\
$51$ & $2\mydot8$ & $0$ & $8\mydot3$ \\
$52$ & $2^3$ & $3$ & $0$ \\
$53$ & $2^4$ & $0$ & $0$ \\
$54$ & $2\mydot4$ & $0$ & $0$ \\
$55$ & $0$ & $3$ & $16\mydot3\mydot5\mydot29$ \\
$56$ & $0$ & $0$ & $2$ \\
$57$ & $2$ & $0$ & $2^2$ \\
$58$ & $2$ & $0$ & $2$ \\
$59$ & $2^2$ & $0$ & $8\mydot9\mydot7\mydot11\mydot31$ \\
$60$ & $4$ & $0$ & $0$ \\
$61$ & $0$ & $0$ & $0$ \\
\hline
$62$ & $2^4$ & $3$ & $0$ \\
$63$ & $2^2\mydot4$ & $0$ & $128\mydot3\mydot5\mydot17$ \\
$64$ & $2^5\mydot4$ & $0$ & $2$ \\
$65$ & $2^7\mydot4$ & $3$ & $2^2$ \\
$66$ & $2^5\mydot8$ & $0$ & $2$ \\
$67$ & $2^3\mydot4$ & $0$ & $8\mydot3$ \\
$68$ & $2^3$ & $3$ & $0$ \\
$69$ & $2^4$ & $0$ & $0$ \\
$70$ & $2^5\mydot4^2$ or $2^6\mydot4$ & $0$ & $0$ \\
$71$ & $2^6\mydot4\mydot8$ or $2^5\mydot4\mydot8$  & $0$ & $16\mydot27\mydot5\mydot7\mydot13\mydot19\mydot37$ \\
$72$ & $2^7$ & $3$ & $2$ \\
$73$ & $2^5$ & $0$ & $2^2$ \\
$74$ & $4^3$ & $3$ & $2$ \\
$75$ & $2$ & $9$ & $8\mydot3$ \\
$76$ & $2^2\mydot4$ & $5$ & $0$ \\
$77$ & $2^5\mydot4$ & $0$ & $0$ \\
$78$ & $2^3\mydot4^2$ & $3$ & $0$ \\
$79$ & $2^6\mydot4$ & $0$ & $32\mydot3\mydot25\mydot11\mydot41$ \\
$80$ & $2^8$ & $0$ & $2$ \\
$81$ & $2^3\mydot4\mydot8$ & $3^2$ & $2^2$ \\
$82$ & $2^5\mydot8$ or $2^4\mydot8$ or  & $3\mydot7$ & $2$ \\
& $2^3\mydot4\mydot8$ \\
$83$ & $2^3\mydot8$ or $2^3\mydot4$ & $5$ & $8\mydot9\mydot49\mydot43$ \\
$84$ & $2^6$ or $2^5$ or $2^4\mydot4$ & $3^2$ & $0$ \\
$85$ & $2^6\mydot4^2$ or $2^5\mydot4^2$ or & $3^2$ & $0$ \\
& $2^4\mydot4^3$ \\
$86$ & $2^5\mydot8^2$ or $2^4\mydot8^2$ or & $3\mydot5$ & $0$ \\
& $2^3\mydot4\mydot8^2$ or $2^2\mydot4\mydot8^2$ \\
$87$ & $2^8$ or $2^7$ or & $0$ & $16\mydot3\mydot5\mydot23$ \\
& $2^6\mydot4$ or $2^5\mydot4$ \\
$88$ & $2^4\mydot4$  & $0$ & $2$ \\
$89$ & $2^3$ & $0$ & $2^2$ \\
$90$ & $2^3\mydot8$ or $2^2\mydot8$ & $3$ & $2$
\end{longtable}

\begin{figure}[b]

\caption{$2$-primary stable homotopy groups
\label{fig:Hatcher}}

\psset{unit=0.43cm}
\renewcommand{\cirrad}{0.15}
\begin{pspicture}(0,-1)(28,11)
\psgrid[unit=2,gridcolor=gridline,subgriddiv=0,gridlabelcolor=white](0.0,0)(14.0,5.0)
\rput(0,-1){0}
\rput(2,-1){2}
\rput(4,-1){4}
\rput(6,-1){6}
\rput(8,-1){8}
\rput(10,-1){10}
\rput(12,-1){12}
\rput(14,-1){14}
\rput(16,-1){16}
\rput(18,-1){18}
\rput(20,-1){20}
\rput(22,-1){22}
\rput(24,-1){24}
\rput(26,-1){26}
\rput(28,-1){28}

\psline[linecolor=vonecolor](0,0)(0,1)
\psline[linecolor=vonecolor](0,0)(1,1)
\psline[linecolor=vonecolor](0,0)(3,1)
\pscircle*[linecolor=vonecolor](0,0){\cirrad}
\psline[linecolor=vonecolor](7,1)(7,2)
\psline[linecolor=vonecolor](7,1)(8,2)
\pscircle*[linecolor=vonecolor](7,1){\cirrad}
\psline[linecolor=vonecolor](9,1)(10,2)
\pscircle*[linecolor=vonecolor](9,1){\cirrad}
\psline[linecolor=vtwocolor](14,6)(15,7)
\psline[linecolor=vtwocolor](14,6)(17,6)
\pscircle*[linecolor=vtwocolor](14,6){\cirrad}
\pscircle*[linecolor=vthreecolor](14,7){\cirrad}
\psline[linecolor=vonecolor](15,1)(15,2)
\psline[linecolor=vonecolor](15,1)(16,2)
\pscircle*[linecolor=vonecolor](15,1){\cirrad}
\psline[linecolor=vthreecolor](16,7)(17,8)
\pscircle*[linecolor=vthreecolor](16,7){\cirrad}
\psline[linecolor=vonecolor](17,1)(18,2)
\pscircle*[linecolor=vonecolor](17,1){\cirrad}
\psline[linecolor=vthreecolor](19,8)(22,8)
\pscircle*[linecolor=vthreecolor](19,8){\cirrad}
\psline[linecolor=vonecolor](23,1)(23,2)
\psline[linecolor=vonecolor](23,1)(24,2)
\pscircle*[linecolor=vonecolor](23,1){\cirrad}
\psline[linecolor=vthreecolor](23,8)(24,9)
\pscircle*[linecolor=vthreecolor](23,8){\cirrad}
\psline[linecolor=vonecolor](25,1)(26,2)
\pscircle*[linecolor=vonecolor](25,1){\cirrad}
\pscircle*[linecolor=vtwocolor](28,5){\cirrad}
\psline[linecolor=vonecolor]{->}(0,3)(0,4.0)
\pscircle*[linecolor=vonecolor](0,3){\cirrad}
\psline[linecolor=vonecolor](0,1)(0,2)
\psline[linecolor=vonecolor](0,1)(3,2)
\pscircle*[linecolor=vonecolor](0,1){\cirrad}
\psline[linecolor=vonecolor](0,2)(0,3)
\psline[linecolor=vonecolor](0,2)(3,3)
\pscircle*[linecolor=vonecolor](0,2){\cirrad}
\psline[linecolor=vonecolor](1,1)(2,2)
\pscircle*[linecolor=vonecolor](1,1){\cirrad}
\psline[linecolor=vonecolor](2,2)(3,3)
\pscircle*[linecolor=vonecolor](2,2){\cirrad}
\psline[linecolor=vonecolor](3,1)(3,2)
\psline[linecolor=vonecolor](3,1)(6,5)
\pscircle*[linecolor=vonecolor](3,1){\cirrad}
\psline[linecolor=vonecolor](3,2)(3,3)
\pscircle*[linecolor=vonecolor](3,2){\cirrad}
\pscircle*[linecolor=vonecolor](3,3){\cirrad}
\psline[linecolor=vtwocolor](6,5)(9,5)
\pscircle*[linecolor=vtwocolor](6,5){\cirrad}
\psline[linecolor=vonecolor](7,2)(7,3)
\pscircle*[linecolor=vonecolor](7,2){\cirrad}
\psline[linecolor=vonecolor](7,3)(7,4)
\pscircle*[linecolor=vonecolor](7,3){\cirrad}
\pscircle*[linecolor=vonecolor](7,4){\cirrad}
\psline[linecolor=vonecolor](8,2)(9,3)
\pscircle*[linecolor=vonecolor](8,2){\cirrad}
\psline[linecolor=vtwocolor](8,4)(9,5)
\pscircle*[linecolor=vtwocolor](8,4){\cirrad}
\pscircle*[linecolor=vonecolor](9,3){\cirrad}
\pscircle*[linecolor=vtwocolor](9,5){\cirrad}
\psline[linecolor=vonecolor](10,2)(11,3)
\pscircle*[linecolor=vonecolor](10,2){\cirrad}
\psline[linecolor=vonecolor](11,1)(11,2)
\pscircle*[linecolor=vonecolor](11,1){\cirrad}
\psline[linecolor=vonecolor](11,2)(11,3)
\pscircle*[linecolor=vonecolor](11,2){\cirrad}
\pscircle*[linecolor=vonecolor](11,3){\cirrad}
\psline[linecolor=vonecolor](15,2)(15,3)
\pscircle*[linecolor=vonecolor](15,2){\cirrad}
\psline[linecolor=vonecolor](15,3)(15,4)
\pscircle*[linecolor=vonecolor](15,3){\cirrad}
\psline[linecolor=vonecolor](15,4)(15,5)
\pscircle*[linecolor=vonecolor](15,4){\cirrad}
\pscircle*[linecolor=vonecolor](15,5){\cirrad}
\pscircle*[linecolor=vtwocolor](15,7){\cirrad}
\psline[linecolor=vonecolor](16,2)(17,3)
\pscircle*[linecolor=vonecolor](16,2){\cirrad}
\pscircle*[linecolor=vonecolor](17,3){\cirrad}
\psline[linecolor=vtwocolor](17,6)(20,6)
\pscircle*[linecolor=vtwocolor](17,6){\cirrad}
\psline[linecolor=vthreecolor](17,8)(18,9)
\pscircle*[linecolor=vthreecolor](17,8){\cirrad}
\psline[linecolor=vonecolor](18,2)(19,3)
\pscircle*[linecolor=vonecolor](18,2){\cirrad}
\psline[linecolor=vthreecolor](18,7)(18,8)
\psline[linecolor=vthreecolor](18,7)(21,7)
\pscircle*[linecolor=vthreecolor](18,7){\cirrad}
\psline[linecolor=vthreecolor](18,8)(18,9)
\pscircle*[linecolor=vthreecolor](18,8){\cirrad}
\pscircle*[linecolor=vthreecolor](18,9){\cirrad}
\psline[linecolor=vonecolor](19,1)(19,2)
\pscircle*[linecolor=vonecolor](19,1){\cirrad}
\psline[linecolor=vonecolor](19,2)(19,3)
\pscircle*[linecolor=vonecolor](19,2){\cirrad}
\pscircle*[linecolor=vonecolor](19,3){\cirrad}
\psline[linecolor=vtwocolor](20,4)(20,5)
\psline[linecolor=vtwocolor](20,4)(21,5)
\psline[linecolor=vtwocolor](20,4)(23,5)
\pscircle*[linecolor=vtwocolor](20,4){\cirrad}
\psline[linecolor=vtwocolor](20,5)(20,6)
\psline[linecolor=vtwocolor](20,5)(23,6)
\pscircle*[linecolor=vtwocolor](20,5){\cirrad}
\psline[linecolor=vtwocolor](20,6)(23,7)
\pscircle*[linecolor=vtwocolor](20,6){\cirrad}
\psline[linecolor=vtwocolor](21,5)(22,6)
\pscircle*[linecolor=vtwocolor](21,5){\cirrad}
\pscircle*[linecolor=vthreecolor](21,7){\cirrad}
\psline[linecolor=vtwocolor](22,6)(23,7)
\pscircle*[linecolor=vtwocolor](22,6){\cirrad}
\pscircle*[linecolor=vthreecolor](22,8){\cirrad}
\psline[linecolor=vonecolor](23,2)(23,3)
\pscircle*[linecolor=vonecolor](23,2){\cirrad}
\psline[linecolor=vonecolor](23,3)(23,4)
\pscircle*[linecolor=vonecolor](23,3){\cirrad}
\pscircle*[linecolor=vonecolor](23,4){\cirrad}
\psline[linecolor=vtwocolor](23,5)(23,6)
\psline[linecolor=vtwocolor](23,5)(26,5)
\pscircle*[linecolor=vtwocolor](23,5){\cirrad}
\psline[linecolor=vtwocolor](23,6)(23,7)
\pscircle*[linecolor=vtwocolor](23,6){\cirrad}
\pscircle*[linecolor=vtwocolor](23,7){\cirrad}
\psline[linecolor=vonecolor](24,2)(25,3)
\pscircle*[linecolor=vonecolor](24,2){\cirrad}
\pscircle*[linecolor=vthreecolor](24,9){\cirrad}
\pscircle*[linecolor=vonecolor](25,3){\cirrad}
\psline[linecolor=vonecolor](26,2)(27,3)
\pscircle*[linecolor=vonecolor](26,2){\cirrad}
\pscircle*[linecolor=vtwocolor](26,5){\cirrad}
\psline[linecolor=vonecolor](27,1)(27,2)
\pscircle*[linecolor=vonecolor](27,1){\cirrad}
\psline[linecolor=vonecolor](27,2)(27,3)
\pscircle*[linecolor=vonecolor](27,2){\cirrad}
\pscircle*[linecolor=vonecolor](27,3){\cirrad}
\end{pspicture}

\end{figure}
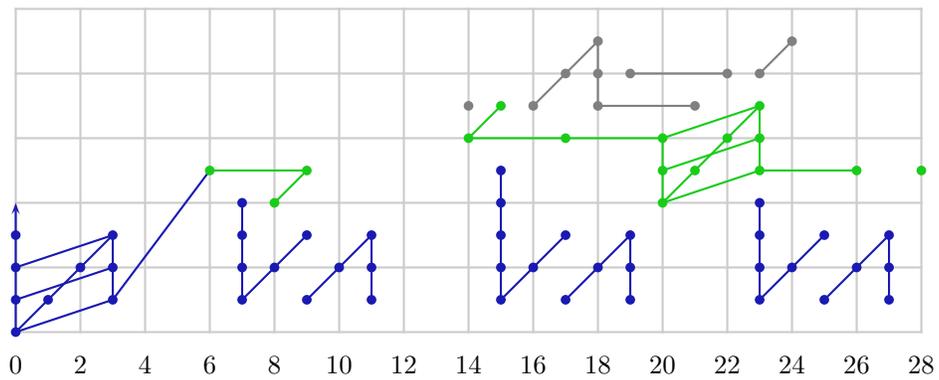

\psset{unit=0.4cm}
\renewcommand{\cirrad}{0.15}
\begin{pspicture}(30,-2)(60,17)
\psgrid[unit=2,gridcolor=gridline,subgriddiv=0,gridlabelcolor=white](15.0,0)(30.0,8.0)
\rput(30,-1){30}
\rput(32,-1){32}
\rput(34,-1){34}
\rput(36,-1){36}
\rput(38,-1){38}
\rput(40,-1){40}
\rput(42,-1){42}
\rput(44,-1){44}
\rput(46,-1){46}
\rput(48,-1){48}
\rput(50,-1){50}
\rput(52,-1){52}
\rput(54,-1){54}
\rput(56,-1){56}
\rput(58,-1){58}
\rput(60,-1){60}

\psline[linecolor=vthreecolor](30,7)(31,8)
\psline[linecolor=vthreecolor](30,7)(33,7)
\pscircle*[linecolor=vthreecolor](30,7){\cirrad}
\psline[linecolor=vonecolor](31,1)(31,2)
\psline[linecolor=vonecolor](31,1)(32,2)
\pscircle*[linecolor=vonecolor](31,1){\cirrad}
\psline[linecolor=vthreecolor](31,9)(34,9)
\pscircle*[linecolor=vthreecolor](31,9){\cirrad}
\psline[linecolor=vtwocolor](32,5)(33,6)
\psline[linecolor=vtwocolor](32,5)(35,5)
\pscircle*[linecolor=vtwocolor](32,5){\cirrad}
\psline[linecolor=vthreecolor](32,8)(35,8)
\pscircle*[linecolor=vthreecolor](32,8){\cirrad}
\psline[linecolor=vthreecolor](32,10)(33,11)
\pscircle*[linecolor=vthreecolor](32,10){\cirrad}
\psline[linecolor=vonecolor](33,1)(34,2)
\pscircle*[linecolor=vonecolor](33,1){\cirrad}
\psline[linecolor=vthreecolor](36,6)(39,6)
\pscircle*[linecolor=vthreecolor](36,6){\cirrad}
\psline[linecolor=vthreecolor](37,9)(40,9)
\pscircle*[linecolor=vthreecolor](37,9){\cirrad}
\psline[linecolor=vonecolor](39,1)(39,2)
\psline[linecolor=vonecolor](39,1)(40,2)
\pscircle*[linecolor=vonecolor](39,1){\cirrad}
\psline[linecolor=vtwocolor](39,5)(40,6)
\psline[linecolor=vtwocolor](39,5)(42,7)
\pscircle*[linecolor=vtwocolor](39,5){\cirrad}
\psline[linecolor=vthreecolor](39,7)(40,8)
\pscircle*[linecolor=vthreecolor](39,7){\cirrad}
\pscircle*[linecolor=vthreecolor](39,10){\cirrad}
\psline[linecolor=vthreecolor](40,7)(41,8)
\pscircle*[linecolor=vthreecolor](40,7){\cirrad}
\psline[linecolor=vthreecolor](40,11)(41,12)
\pscircle*[linecolor=vthreecolor](40,11){\cirrad}
\psline[linecolor=vonecolor](41,1)(42,2)
\pscircle*[linecolor=vonecolor](41,1){\cirrad}
\psline[linecolor=vthreecolor](44,14)(44,15)
\psline[linecolor=vthreecolor](44,14)(45,15)
\psline[linecolor=vthreecolor](44,14)(47,14)
\pscircle*[linecolor=vthreecolor](44,14){\cirrad}
\psline[linecolor=vtwocolor](45,6)(46,7)
\psline[linecolor=vtwocolor](45,6)(48,7)
\pscircle*[linecolor=vtwocolor](45,6){\cirrad}
\psline[linecolor=vthreecolor](45,12)(46,13)
\psline[linecolor=vthreecolor](45,12)(48,12)
\pscircle*[linecolor=vthreecolor](45,12){\cirrad}
\psline[linecolor=vonecolor](47,1)(47,2)
\psline[linecolor=vonecolor](47,1)(48,2)
\pscircle*[linecolor=vonecolor](47,1){\cirrad}
\psline[linecolor=vthreecolor](47,13)(48,14)
\pscircle*[linecolor=vthreecolor](47,13){\cirrad}
\psline[linecolor=vonecolor](49,1)(50,2)
\pscircle*[linecolor=vonecolor](49,1){\cirrad}
\psline[linecolor=vthreecolor](50,7)(53,7)
\pscircle*[linecolor=vthreecolor](50,7){\cirrad}
\psline[linecolor=vtwocolor](52,4)(53,5)
\pscircle*[linecolor=vtwocolor](52,4){\cirrad}
\pscircle*[linecolor=vthreecolor](52,8){\cirrad}
\psline[linecolor=vtwocolor](54,5)(54,6)
\psline[linecolor=vtwocolor](54,5)(57,5)
\pscircle*[linecolor=vtwocolor](54,5){\cirrad}
\psline[linecolor=vonecolor](55,1)(55,2)
\psline[linecolor=vonecolor](55,1)(56,2)
\pscircle*[linecolor=vonecolor](55,1){\cirrad}
\psline[linecolor=vonecolor](57,1)(58,2)
\pscircle*[linecolor=vonecolor](57,1){\cirrad}
\pscircle*[linecolor=vthreecolor](58,6){\cirrad}
\psline[linecolor=vtwocolor](59,4)(60,5)
\pscircle*[linecolor=vtwocolor](59,4){\cirrad}
\pscircle*[linecolor=vthreecolor](59,6){\cirrad}
\psline[linecolor=vonecolor](31,2)(31,3)
\pscircle*[linecolor=vonecolor](31,2){\cirrad}
\psline[linecolor=vonecolor](31,3)(31,4)
\pscircle*[linecolor=vonecolor](31,3){\cirrad}
\psline[linecolor=vonecolor](31,4)(31,5)
\pscircle*[linecolor=vonecolor](31,4){\cirrad}
\psline[linecolor=vonecolor](31,5)(31,6)
\pscircle*[linecolor=vonecolor](31,5){\cirrad}
\pscircle*[linecolor=vonecolor](31,6){\cirrad}
\pscircle*[linecolor=vthreecolor](31,8){\cirrad}
\psline[linecolor=vonecolor](32,2)(33,3)
\pscircle*[linecolor=vonecolor](32,2){\cirrad}
\pscircle*[linecolor=vonecolor](33,3){\cirrad}
\pscircle*[linecolor=vtwocolor](33,6){\cirrad}
\psline[linecolor=vthreecolor](33,11)(34,12)
\pscircle*[linecolor=vthreecolor](33,11){\cirrad}
\pscircle*[linecolor=vthreecolor](33,7){\cirrad}
\psline[linecolor=vonecolor](34,2)(35,3)
\pscircle*[linecolor=vonecolor](34,2){\cirrad}
\psline[linecolor=vtwocolor](34,4)(35,5)
\pscircle*[linecolor=vtwocolor](34,4){\cirrad}
\psline[linecolor=vthreecolor](34,11)(34,12)
\pscircle*[linecolor=vthreecolor](34,11){\cirrad}
\pscircle*[linecolor=vthreecolor](34,12){\cirrad}
\pscircle*[linecolor=vthreecolor](34,9){\cirrad}
\psline[linecolor=vonecolor](35,1)(35,2)
\pscircle*[linecolor=vonecolor](35,1){\cirrad}
\psline[linecolor=vonecolor](35,2)(35,3)
\pscircle*[linecolor=vonecolor](35,2){\cirrad}
\pscircle*[linecolor=vonecolor](35,3){\cirrad}
\pscircle*[linecolor=vtwocolor](35,5){\cirrad}
\psline[linecolor=vthreecolor](35,8)(38,8)
\pscircle*[linecolor=vthreecolor](35,8){\cirrad}
\psline[linecolor=vthreecolor](37,7)(38,8)
\pscircle*[linecolor=vthreecolor](37,7){\cirrad}
\pscircle*[linecolor=vthreecolor](38,8){\cirrad}
\psline[linecolor=vthreecolor](38,4)(38,5)
\psline[linecolor=vthreecolor](38,4)(39,6)
\pscircle*[linecolor=vthreecolor](38,4){\cirrad}
\pscircle*[linecolor=vthreecolor](38,5){\cirrad}
\psline[linecolor=vonecolor](39,2)(39,3)
\pscircle*[linecolor=vonecolor](39,2){\cirrad}
\psline[linecolor=vonecolor](39,3)(39,4)
\pscircle*[linecolor=vonecolor](39,3){\cirrad}
\pscircle*[linecolor=vonecolor](39,4){\cirrad}
\psline[linecolor=vthreecolor](39,8)(40,9)
\pscircle*[linecolor=vthreecolor](39,8){\cirrad}
\pscircle*[linecolor=vthreecolor](39,6){\cirrad}
\psline[linecolor=vonecolor](40,2)(41,3)
\pscircle*[linecolor=vonecolor](40,2){\cirrad}
\psline[linecolor=vtwocolor](40,5)(40,6)
\psline[linecolor=vtwocolor](40,5)(41,6)
\pscircle*[linecolor=vtwocolor](40,5){\cirrad}
\pscircle*[linecolor=vtwocolor](40,6){\cirrad}
\pscircle*[linecolor=vthreecolor](40,9){\cirrad}
\pscircle*[linecolor=vthreecolor](40,8){\cirrad}
\pscircle*[linecolor=vonecolor](41,3){\cirrad}
\psline[linecolor=vtwocolor](41,6)(42,7)
\pscircle*[linecolor=vtwocolor](41,6){\cirrad}
\pscircle*[linecolor=vthreecolor](41,8){\cirrad}
\psline[linecolor=vthreecolor](41,12)(42,13)
\pscircle*[linecolor=vthreecolor](41,12){\cirrad}
\psline[linecolor=vonecolor](42,2)(43,3)
\pscircle*[linecolor=vonecolor](42,2){\cirrad}
\pscircle*[linecolor=vtwocolor](42,7){\cirrad}
\psline[linecolor=vthreecolor](42,11)(42,12)
\psline[linecolor=vthreecolor](42,11)(45,11)
\pscircle*[linecolor=vthreecolor](42,11){\cirrad}
\psline[linecolor=vthreecolor](42,12)(42,13)
\pscircle*[linecolor=vthreecolor](42,12){\cirrad}
\pscircle*[linecolor=vthreecolor](42,13){\cirrad}
\psline[linecolor=vonecolor](43,1)(43,2)
\pscircle*[linecolor=vonecolor](43,1){\cirrad}
\psline[linecolor=vonecolor](43,2)(43,3)
\pscircle*[linecolor=vonecolor](43,2){\cirrad}
\pscircle*[linecolor=vonecolor](43,3){\cirrad}
\psline[linecolor=vthreecolor](44,15)(44,16)
\pscircle*[linecolor=vthreecolor](44,15){\cirrad}
\pscircle*[linecolor=vthreecolor](44,16){\cirrad}
\psline[linecolor=vthreecolor](45,15)(46,16)
\pscircle*[linecolor=vthreecolor](45,15){\cirrad}
\psline[linecolor=vthreecolor](45,8)(45,9)
\psline[linecolor=vthreecolor](45,8)(46,9)
\psline[linecolor=vthreecolor](45,8)(48,9)
\pscircle*[linecolor=vthreecolor](45,8){\cirrad}
\psline[linecolor=vthreecolor](45,9)(45,10)
\psline[linecolor=vthreecolor](45,9)(48,10)
\pscircle*[linecolor=vthreecolor](45,9){\cirrad}
\psline[linecolor=vthreecolor](45,10)(45,11)
\pscircle*[linecolor=vthreecolor](45,10){\cirrad}
\pscircle*[linecolor=vthreecolor](45,11){\cirrad}
\psline[linecolor=vtwocolor](46,7)(47,8)
\pscircle*[linecolor=vtwocolor](46,7){\cirrad}
\pscircle*[linecolor=vthreecolor](46,16){\cirrad}
\pscircle*[linecolor=vthreecolor](46,13){\cirrad}
\psline[linecolor=vthreecolor](46,9)(47,10)
\pscircle*[linecolor=vthreecolor](46,9){\cirrad}
\psline[linecolor=vonecolor](47,2)(47,3)
\pscircle*[linecolor=vonecolor](47,2){\cirrad}
\psline[linecolor=vonecolor](47,3)(47,4)
\pscircle*[linecolor=vonecolor](47,3){\cirrad}
\psline[linecolor=vonecolor](47,4)(47,5)
\pscircle*[linecolor=vonecolor](47,4){\cirrad}
\pscircle*[linecolor=vonecolor](47,5){\cirrad}
\psline[linecolor=vtwocolor](47,6)(47,8)
\psline[linecolor=vtwocolor](47,6)(48,7)
\pscircle*[linecolor=vtwocolor](47,6){\cirrad}
\pscircle*[linecolor=vtwocolor](47,8){\cirrad}
\pscircle*[linecolor=vthreecolor](47,14){\cirrad}
\pscircle*[linecolor=vthreecolor](47,10){\cirrad}
\psline[linecolor=vonecolor](48,2)(49,3)
\pscircle*[linecolor=vonecolor](48,2){\cirrad}
\pscircle*[linecolor=vtwocolor](48,7){\cirrad}
\psline[linecolor=vthreecolor](48,12)(51,12)
\pscircle*[linecolor=vthreecolor](48,12){\cirrad}
\psline[linecolor=vthreecolor](48,9)(48,10)
\psline[linecolor=vthreecolor](48,9)(51,9)
\pscircle*[linecolor=vthreecolor](48,9){\cirrad}
\pscircle*[linecolor=vthreecolor](48,10){\cirrad}
\pscircle*[linecolor=vthreecolor](48,14){\cirrad}
\pscircle*[linecolor=vonecolor](49,3){\cirrad}
\psline[linecolor=vonecolor](50,2)(51,3)
\pscircle*[linecolor=vonecolor](50,2){\cirrad}
\psline[linecolor=vthreecolor](50,11)(51,12)
\psline[linecolor=vthreecolor](50,11)(53,12)
\pscircle*[linecolor=vthreecolor](50,11){\cirrad}
\psline[linecolor=vonecolor](51,1)(51,2)
\pscircle*[linecolor=vonecolor](51,1){\cirrad}
\psline[linecolor=vonecolor](51,2)(51,3)
\pscircle*[linecolor=vonecolor](51,2){\cirrad}
\pscircle*[linecolor=vonecolor](51,3){\cirrad}
\psline[linecolor=vthreecolor](51,10)(51,11)
\psline[linecolor=vthreecolor](51,10)(52,11)
\pscircle*[linecolor=vthreecolor](51,10){\cirrad}
\psline[linecolor=vthreecolor](51,11)(51,12)
\pscircle*[linecolor=vthreecolor](51,11){\cirrad}
\pscircle*[linecolor=vthreecolor](51,12){\cirrad}
\psline[linecolor=vthreecolor](51,9)(54,9)
\pscircle*[linecolor=vthreecolor](51,9){\cirrad}
\pscircle*[linecolor=vthreecolor](52,11){\cirrad}
\pscircle*[linecolor=vtwocolor](53,5){\cirrad}
\pscircle*[linecolor=vthreecolor](53,7){\cirrad}
\pscircle*[linecolor=vthreecolor](53,12){\cirrad}
\psline[linecolor=vthreecolor](53,8)(54,9)
\pscircle*[linecolor=vthreecolor](53,8){\cirrad}
\pscircle*[linecolor=vtwocolor](54,6){\cirrad}
\pscircle*[linecolor=vthreecolor](54,9){\cirrad}
\psline[linecolor=vonecolor](55,2)(55,3)
\pscircle*[linecolor=vonecolor](55,2){\cirrad}
\psline[linecolor=vonecolor](55,3)(55,4)
\pscircle*[linecolor=vonecolor](55,3){\cirrad}
\pscircle*[linecolor=vonecolor](55,4){\cirrad}
\psline[linecolor=vonecolor](56,2)(57,3)
\pscircle*[linecolor=vonecolor](56,2){\cirrad}
\pscircle*[linecolor=vonecolor](57,3){\cirrad}
\psline[linecolor=vtwocolor](57,5)(60,5)
\pscircle*[linecolor=vtwocolor](57,5){\cirrad}
\psline[linecolor=vonecolor](58,2)(59,3)
\pscircle*[linecolor=vonecolor](58,2){\cirrad}
\psline[linecolor=vonecolor](59,1)(59,2)
\pscircle*[linecolor=vonecolor](59,1){\cirrad}
\psline[linecolor=vonecolor](59,2)(59,3)
\pscircle*[linecolor=vonecolor](59,2){\cirrad}
\pscircle*[linecolor=vonecolor](59,3){\cirrad}
\psline[linecolor=vtwocolor](60,4)(60,5)
\pscircle*[linecolor=vtwocolor](60,4){\cirrad}
\pscircle*[linecolor=vtwocolor](60,5){\cirrad}
\end{pspicture}

\psset{unit=0.43cm}
\renewcommand{\cirrad}{0.15}
\begin{pspicture}(62,-2)(90,23)
\psgrid[unit=2,gridcolor=gridline,subgriddiv=0,gridlabelcolor=white](31.0,0)(45.0,11.0)
\rput(62,-1){62}
\rput(64,-1){64}
\rput(66,-1){66}
\rput(68,-1){68}
\rput(70,-1){70}
\rput(72,-1){72}
\rput(74,-1){74}
\rput(76,-1){76}
\rput(78,-1){78}
\rput(80,-1){80}
\rput(82,-1){82}
\rput(84,-1){84}
\rput(86,-1){86}
\rput(88,-1){88}
\rput(90,-1){90}

\pscircle*[linecolor=vthreecolor](62,8){\cirrad}
\psline[linecolor=vthreecolor](62,12)(65,12)
\pscircle*[linecolor=vthreecolor](62,12){\cirrad}
\psline[linecolor=vthreecolor](62,13)(63,14)
\psline[linecolor=vthreecolor](62,13)(65,13)
\pscircle*[linecolor=vthreecolor](62,13){\cirrad}
\psline[linecolor=vthreecolor](62,19)(63,20)
\psline[linecolor=vthreecolor](62,19)(65,20)
\pscircle*[linecolor=vthreecolor](62,19){\cirrad}
\psline[linecolor=vonecolor](63,1)(63,2)
\psline[linecolor=vonecolor](63,1)(64,2)
\pscircle*[linecolor=vonecolor](63,1){\cirrad}
\psline[linecolor=vthreecolor](64,9)(65,10)
\pscircle*[linecolor=vthreecolor](64,9){\cirrad}
\psline[linecolor=vthreecolor](64,8)(67,8)
\pscircle*[linecolor=vthreecolor](64,8){\cirrad}
\psline[linecolor=vonecolor](65,1)(66,2)
\pscircle*[linecolor=vonecolor](65,1){\cirrad}
\psline[linecolor=vtwocolor](65,5)(66,6)
\pscircle*[linecolor=vtwocolor](65,5){\cirrad}
\psline[linecolor=vtwocolor](65,7)(68,7)
\pscircle*[linecolor=vtwocolor](65,7){\cirrad}
\pscircle*[linecolor=vthreecolor](66,12){\cirrad}
\psline[linecolor=vthreecolor](67,15)(70,15)
\pscircle*[linecolor=vthreecolor](67,15){\cirrad}
\pscircle*[linecolor=vthreecolor](70,10){\cirrad}
\psline[linecolor=vthreecolor](70,14)(71,15)
\psline[linecolor=vthreecolor](70,14)(73,14)
\pscircle*[linecolor=vthreecolor](70,14){\cirrad}
\psline[linecolor=vonecolor](71,1)(71,2)
\psline[linecolor=vonecolor](71,1)(72,2)
\pscircle*[linecolor=vonecolor](71,1){\cirrad}
\psline[linecolor=vthreecolor](71,8)(72,9)
\pscircle*[linecolor=vthreecolor](71,8){\cirrad}
\pscircle*[linecolor=vthreecolor](71,13){\cirrad}
\psline[linecolor=vzerocolor](71,16)(72,17)
\pscircle*[linecolor=vzerocolor](71,16){\cirrad}
\psline[linecolor=vthreecolor](72,10)(73,11)
\pscircle*[linecolor=vthreecolor](72,10){\cirrad}
\psline[linecolor=vonecolor](73,1)(74,2)
\pscircle*[linecolor=vonecolor](73,1){\cirrad}
\pscircle*[linecolor=vthreecolor](73,8){\cirrad}
\psline[linecolor=vthreecolor](74,8)(74,9)
\pscircle*[linecolor=vthreecolor](74,8){\cirrad}
\psline[linecolor=vthreecolor](75,8)(76,9)
\pscircle*[linecolor=vthreecolor](75,8){\cirrad}
\psline[linecolor=vthreecolor](76,7)(76,8)
\psline[linecolor=vthreecolor](76,7)(79,7)
\pscircle*[linecolor=vthreecolor](76,7){\cirrad}
\psline[linecolor=vthreecolor](76,13)(77,14)
\psline[linecolor=vthreecolor](76,13)(79,13)
\pscircle*[linecolor=vthreecolor](76,13){\cirrad}
\psline[linecolor=vthreecolor](77,8)(80,8)
\pscircle*[linecolor=vthreecolor](77,8){\cirrad}
\pscircle*[linecolor=vthreecolor](77,9){\cirrad}
\psline[linecolor=vthreecolor](77,16)(78,17)
\psline[linecolor=vthreecolor](77,16)(80,16)
\pscircle*[linecolor=vthreecolor](77,16){\cirrad}
\psline[linecolor=vthreecolor](77,18)(78,19)
\psline[linecolor=vthreecolor](77,18)(80,18)
\pscircle*[linecolor=vthreecolor](77,18){\cirrad}
\psline[linecolor=vthreecolor](78,10)(78,11)
\psline[linecolor=vthreecolor](78,10)(81,10)
\pscircle*[linecolor=vthreecolor](78,10){\cirrad}
\psline[linecolor=vonecolor](79,1)(79,2)
\psline[linecolor=vonecolor](79,1)(80,2)
\pscircle*[linecolor=vonecolor](79,1){\cirrad}
\pscircle*[linecolor=vthreecolor](79,15){\cirrad}
\pscircle*[linecolor=vthreecolor](79,9){\cirrad}
\psline[linecolor=vzerocolor](79,17)(82,17)
\pscircle*[linecolor=vthreecolor](79,17){\cirrad}
\psline[linecolor=vthreecolor](79,19)(80,20)
\pscircle*[linecolor=vthreecolor](79,19){\cirrad}
\psline[linecolor=vthreecolor](79,12)(80,13)
\pscircle*[linecolor=vthreecolor](79,12){\cirrad}
\pscircle*[linecolor=vtwocolor](80,5){\cirrad}
\psline[linecolor=vonecolor](81,1)(82,2)
\pscircle*[linecolor=vonecolor](81,1){\cirrad}
\psline[linecolor=vthreecolor](81,19)(84,19)
\pscircle*[linecolor=vthreecolor](81,19){\cirrad}
\pscircle*[linecolor=vthreecolor](81,15){\cirrad}
\psline[linecolor=vthreecolor](82,13)(85,13)
\pscircle*[linecolor=vthreecolor](82,13){\cirrad}
\pscircle*[linecolor=vthreecolor](82,15){\cirrad}
\pscircle*[linecolor=vthreecolor](82,11){\cirrad}
\psline[linecolor=vzerocolor](83,21)(84,22)
\psline[linecolor=vzerocolor](83,21)(86,21)
\pscircle*[linecolor=vthreecolor](83,21){\cirrad}
\psline[linecolor=vthreecolor](84,14)(87,14)
\pscircle*[linecolor=vthreecolor](84,14){\cirrad}
\pscircle*[linecolor=vtwocolor](85,5){\cirrad}
\psline[linecolor=vthreecolor](85,18)(85,19)
\psline[linecolor=vthreecolor](85,18)(86,20)
\psline[linecolor=vthreecolor](85,18)(88,18)
\pscircle*[linecolor=vthreecolor](85,18){\cirrad}
\psline[linecolor=vzerocolor](85,7)(86,8)
\pscircle*[linecolor=vzerocolor](85,7){\cirrad}
\psline[linecolor=vthreecolor](86,5)(87,6)
\pscircle*[linecolor=vthreecolor](86,5){\cirrad}
\psline[linecolor=vonecolor](87,1)(87,2)
\psline[linecolor=vonecolor](87,1)(88,2)
\pscircle*[linecolor=vonecolor](87,1){\cirrad}
\pscircle*[linecolor=vthreecolor](87,16){\cirrad}
\pscircle*[linecolor=vthreecolor](87,17){\cirrad}
\psline[linecolor=vzerocolor](87,7)(88,8)
\pscircle*[linecolor=vzerocolor](87,7){\cirrad}
\pscircle*[linecolor=vthreecolor](87,10){\cirrad}
\psline[linecolor=vthreecolor](87,9)(88,10)
\pscircle*[linecolor=vthreecolor](87,9){\cirrad}
\psline[linecolor=vzerocolor](87,11)(88,12)
\pscircle*[linecolor=vzerocolor](87,11){\cirrad}
\psline[linecolor=vthreecolor](88,6)(89,7)
\pscircle*[linecolor=vthreecolor](88,6){\cirrad}
\psline[linecolor=vonecolor](89,1)(90,2)
\pscircle*[linecolor=vonecolor](89,1){\cirrad}
\pscircle*[linecolor=vtwocolor](90,5){\cirrad}
\pscircle*[linecolor=vzerocolor](90,10){\cirrad}
\pscircle*[linecolor=vthreecolor](90,12){\cirrad}
\psline[linecolor=vonecolor](63,2)(63,3)
\pscircle*[linecolor=vonecolor](63,2){\cirrad}
\psline[linecolor=vonecolor](63,3)(63,4)
\pscircle*[linecolor=vonecolor](63,3){\cirrad}
\psline[linecolor=vonecolor](63,4)(63,5)
\pscircle*[linecolor=vonecolor](63,4){\cirrad}
\psline[linecolor=vonecolor](63,5)(63,6)
\pscircle*[linecolor=vonecolor](63,5){\cirrad}
\psline[linecolor=vonecolor](63,6)(63,7)
\pscircle*[linecolor=vonecolor](63,6){\cirrad}
\pscircle*[linecolor=vonecolor](63,7){\cirrad}
\psline[linecolor=vthreecolor](63,17)(64,19)
\psline[linecolor=vthreecolor](63,17)(66,17)
\pscircle*[linecolor=vthreecolor](63,17){\cirrad}
\psline[linecolor=vthreecolor](63,19)(63,20)
\psline[linecolor=vthreecolor](63,19)(64,20)
\psline[linecolor=vthreecolor](63,19)(66,20)
\pscircle*[linecolor=vthreecolor](63,19){\cirrad}
\pscircle*[linecolor=vthreecolor](63,20){\cirrad}
\psline[linecolor=vthreecolor](63,14)(64,15)
\pscircle*[linecolor=vthreecolor](63,14){\cirrad}
\psline[linecolor=vonecolor](64,2)(65,3)
\pscircle*[linecolor=vonecolor](64,2){\cirrad}
\psline[linecolor=vthreecolor](64,19)(65,20)
\pscircle*[linecolor=vthreecolor](64,19){\cirrad}
\psline[linecolor=vthreecolor](64,20)(65,21)
\pscircle*[linecolor=vthreecolor](64,20){\cirrad}
\psline[linecolor=vthreecolor](64,18)(67,18)
\pscircle*[linecolor=vthreecolor](64,18){\cirrad}
\psline[linecolor=vthreecolor](64,14)(64,15)
\psline[linecolor=vthreecolor](64,14)(65,15)
\pscircle*[linecolor=vthreecolor](64,14){\cirrad}
\pscircle*[linecolor=vthreecolor](64,15){\cirrad}
\pscircle*[linecolor=vonecolor](65,3){\cirrad}
\psline[linecolor=vthreecolor](65,19)(65,20)
\psline[linecolor=vthreecolor](65,19)(66,20)
\pscircle*[linecolor=vthreecolor](65,19){\cirrad}
\pscircle*[linecolor=vthreecolor](65,20){\cirrad}
\pscircle*[linecolor=vthreecolor](65,21){\cirrad}
\psline[linecolor=vthreecolor](65,10)(66,11)
\pscircle*[linecolor=vthreecolor](65,10){\cirrad}
\pscircle*[linecolor=vthreecolor](65,12){\cirrad}
\psline[linecolor=vthreecolor](65,13)(68,14)
\pscircle*[linecolor=vthreecolor](65,13){\cirrad}
\psline[linecolor=vthreecolor](65,15)(66,16)
\pscircle*[linecolor=vthreecolor](65,15){\cirrad}
\psline[linecolor=vonecolor](66,2)(67,3)
\pscircle*[linecolor=vonecolor](66,2){\cirrad}
\pscircle*[linecolor=vtwocolor](66,6){\cirrad}
\psline[linecolor=vthreecolor,linestyle=dashed](66,17)(69,17)
\pscircle*[linecolor=vthreecolor](66,17){\cirrad}
\pscircle*[linecolor=vthreecolor](66,20){\cirrad}
\psline[linecolor=vthreecolor](66,10)(66,11)
\pscircle*[linecolor=vthreecolor](66,10){\cirrad}
\pscircle*[linecolor=vthreecolor](66,11){\cirrad}
\psline[linecolor=vthreecolor](66,9)(66,10)
\psline[linecolor=vthreecolor](66,9)(69,9)
\pscircle*[linecolor=vthreecolor](66,9){\cirrad}
\pscircle*[linecolor=vthreecolor](66,16){\cirrad}
\psline[linecolor=vonecolor](67,1)(67,2)
\pscircle*[linecolor=vonecolor](67,1){\cirrad}
\psline[linecolor=vonecolor](67,2)(67,3)
\pscircle*[linecolor=vonecolor](67,2){\cirrad}
\pscircle*[linecolor=vonecolor](67,3){\cirrad}
\psline[linecolor=vzerocolor](67,18)(70,18)
\pscircle*[linecolor=vthreecolor](67,18){\cirrad}
\pscircle*[linecolor=vthreecolor](67,8){\cirrad}
\psline[linecolor=vthreecolor](67,12)(67,13)
\psline[linecolor=vthreecolor](67,12)(68,14)
\psline[linecolor=vthreecolor](67,12)(70,12)
\pscircle*[linecolor=vthreecolor](67,12){\cirrad}
\psline[linecolor=vthreecolor](67,13)(70,13)
\pscircle*[linecolor=vthreecolor](67,13){\cirrad}
\psline[linecolor=vtwocolor](68,7)(71,7)
\pscircle*[linecolor=vtwocolor](68,7){\cirrad}
\psline[linecolor=vthreecolor](68,16)(69,17)
\pscircle*[linecolor=vthreecolor](68,16){\cirrad}
\pscircle*[linecolor=vthreecolor](68,14){\cirrad}
\psline[linecolor=vzerocolor](69,17)(70,18)
\pscircle*[linecolor=vthreecolor](69,17){\cirrad}
\pscircle*[linecolor=vthreecolor](69,9){\cirrad}
\psline[linecolor=vthreecolor](69,8)(70,9)
\pscircle*[linecolor=vthreecolor](69,8){\cirrad}
\psline[linecolor=vthreecolor](69,11)(70,13)
\psline[linecolor=vthreecolor](69,11)(72,12)
\pscircle*[linecolor=vthreecolor](69,11){\cirrad}
\pscircle*[linecolor=vzerocolor](70,18){\cirrad}
\psline[linecolor=vtwocolor](70,6)(71,7)
\psline[linecolor=vtwocolor](70,6)(73,7)
\pscircle*[linecolor=vtwocolor](70,6){\cirrad}
\psline[linecolor=vzerocolor](70,17)(70,18)
\psline[linecolor=vthreecolor](70,17)(71,18)
\pscircle*[linecolor=vthreecolor](70,17){\cirrad}
\pscircle*[linecolor=vthreecolor](70,15){\cirrad}
\psline[linecolor=vthreecolor](70,12)(70,13)
\psline[linecolor=vthreecolor](70,12)(73,13)
\pscircle*[linecolor=vthreecolor](70,12){\cirrad}
\pscircle*[linecolor=vthreecolor](70,13){\cirrad}
\psline[linecolor=vthreecolor](70,9)(71,10)
\pscircle*[linecolor=vthreecolor](70,9){\cirrad}
\psline[linecolor=vonecolor](71,2)(71,3)
\pscircle*[linecolor=vonecolor](71,2){\cirrad}
\psline[linecolor=vonecolor](71,3)(71,4)
\pscircle*[linecolor=vonecolor](71,3){\cirrad}
\pscircle*[linecolor=vonecolor](71,4){\cirrad}
\psline[linecolor=vtwocolor](71,5)(71,6)
\psline[linecolor=vtwocolor](71,5)(72,6)
\psline[linecolor=vtwocolor](71,5)(74,6)
\pscircle*[linecolor=vtwocolor](71,5){\cirrad}
\psline[linecolor=vtwocolor](71,6)(71,7)
\psline[linecolor=vtwocolor](71,6)(74,7)
\pscircle*[linecolor=vtwocolor](71,6){\cirrad}
\pscircle*[linecolor=vtwocolor](71,7){\cirrad}
\pscircle*[linecolor=vthreecolor](71,18){\cirrad}
\pscircle*[linecolor=vthreecolor](71,15){\cirrad}
\psline[linecolor=vthreecolor](71,11)(72,12)
\pscircle*[linecolor=vthreecolor](71,11){\cirrad}
\psline[linecolor=vthreecolor](71,9)(71,10)
\psline[linecolor=vthreecolor](71,9)(72,11)
\pscircle*[linecolor=vthreecolor](71,9){\cirrad}
\pscircle*[linecolor=vthreecolor](71,10){\cirrad}
\psline[linecolor=vonecolor](72,2)(73,3)
\pscircle*[linecolor=vonecolor](72,2){\cirrad}
\psline[linecolor=vtwocolor](72,6)(73,7)
\pscircle*[linecolor=vtwocolor](72,6){\cirrad}
\pscircle*[linecolor=vthreecolor](72,9){\cirrad}
\psline[linecolor=vthreecolor](72,13)(73,14)
\pscircle*[linecolor=vthreecolor](72,13){\cirrad}
\pscircle*[linecolor=vthreecolor](72,17){\cirrad}
\pscircle*[linecolor=vthreecolor](72,12){\cirrad}
\psline[linecolor=vthreecolor](72,11)(73,13)
\pscircle*[linecolor=vthreecolor](72,11){\cirrad}
\pscircle*[linecolor=vonecolor](73,3){\cirrad}
\pscircle*[linecolor=vtwocolor](73,7){\cirrad}
\pscircle*[linecolor=vthreecolor](73,14){\cirrad}
\pscircle*[linecolor=vthreecolor](73,13){\cirrad}
\psline[linecolor=vthreecolor](73,11)(74,12)
\pscircle*[linecolor=vthreecolor](73,11){\cirrad}
\psline[linecolor=vonecolor](74,2)(75,3)
\pscircle*[linecolor=vonecolor](74,2){\cirrad}
\psline[linecolor=vtwocolor](74,6)(74,7)
\pscircle*[linecolor=vtwocolor](74,6){\cirrad}
\pscircle*[linecolor=vtwocolor](74,7){\cirrad}
\pscircle*[linecolor=vthreecolor](74,9){\cirrad}
\psline[linecolor=vthreecolor](74,11)(74,12)
\pscircle*[linecolor=vthreecolor](74,11){\cirrad}
\pscircle*[linecolor=vthreecolor](74,12){\cirrad}
\psline[linecolor=vonecolor](75,1)(75,2)
\pscircle*[linecolor=vonecolor](75,1){\cirrad}
\psline[linecolor=vonecolor](75,2)(75,3)
\pscircle*[linecolor=vonecolor](75,2){\cirrad}
\pscircle*[linecolor=vonecolor](75,3){\cirrad}
\pscircle*[linecolor=vthreecolor](76,9){\cirrad}
\pscircle*[linecolor=vthreecolor](76,8){\cirrad}
\psline[linecolor=vthreecolor](77,11)(77,12)
\psline[linecolor=vthreecolor](77,11)(78,12)
\pscircle*[linecolor=vthreecolor](77,11){\cirrad}
\pscircle*[linecolor=vthreecolor](77,12){\cirrad}
\psline[linecolor=vthreecolor](77,14)(78,15)
\pscircle*[linecolor=vthreecolor](77,14){\cirrad}
\psline[linecolor=vthreecolor](78,11)(81,11)
\pscircle*[linecolor=vthreecolor](78,11){\cirrad}
\pscircle*[linecolor=vthreecolor](78,19){\cirrad}
\pscircle*[linecolor=vthreecolor](78,17){\cirrad}
\psline[linecolor=vthreecolor](78,12)(79,13)
\pscircle*[linecolor=vthreecolor](78,12){\cirrad}
\psline[linecolor=vthreecolor](78,14)(78,15)
\pscircle*[linecolor=vthreecolor](78,14){\cirrad}
\pscircle*[linecolor=vthreecolor](78,15){\cirrad}
\psline[linecolor=vonecolor](79,2)(79,3)
\pscircle*[linecolor=vonecolor](79,2){\cirrad}
\psline[linecolor=vonecolor](79,3)(79,4)
\pscircle*[linecolor=vonecolor](79,3){\cirrad}
\psline[linecolor=vonecolor](79,4)(79,5)
\pscircle*[linecolor=vonecolor](79,4){\cirrad}
\pscircle*[linecolor=vonecolor](79,5){\cirrad}
\psline[linecolor=vthreecolor](79,6)(79,7)
\psline[linecolor=vthreecolor](79,6)(82,6)
\pscircle*[linecolor=vthreecolor](79,6){\cirrad}
\pscircle*[linecolor=vthreecolor](79,7){\cirrad}
\pscircle*[linecolor=vthreecolor](79,13){\cirrad}
\psline[linecolor=vonecolor](80,2)(81,3)
\pscircle*[linecolor=vonecolor](80,2){\cirrad}
\psline[linecolor=vthreecolor](80,7)(81,8)
\pscircle*[linecolor=vthreecolor](80,7){\cirrad}
\psline[linecolor=vthreecolor](80,9)(81,11)
\psline[linecolor=vthreecolor](80,9)(83,10)
\pscircle*[linecolor=vthreecolor](80,9){\cirrad}
\psline[linecolor=vthreecolor](80,13)(81,14)
\pscircle*[linecolor=vthreecolor](80,13){\cirrad}
\psline[linecolor=vthreecolor](80,18)(83,18)
\pscircle*[linecolor=vthreecolor](80,18){\cirrad}
\pscircle*[linecolor=vthreecolor](80,8){\cirrad}
\psline[linecolor=vthreecolor](80,16)(83,16)
\pscircle*[linecolor=vthreecolor](80,16){\cirrad}
\pscircle*[linecolor=vthreecolor](80,20){\cirrad}
\pscircle*[linecolor=vonecolor](81,3){\cirrad}
\psline[linecolor=vthreecolor](81,8)(82,9)
\pscircle*[linecolor=vthreecolor](81,8){\cirrad}
\psline[linecolor=vthreecolor](81,10)(81,11)
\pscircle*[linecolor=vthreecolor](81,10){\cirrad}
\pscircle*[linecolor=vthreecolor](81,11){\cirrad}
\psline[linecolor=vthreecolor](81,12)(81,13)
\psline[linecolor=vthreecolor](81,12)(84,12)
\pscircle*[linecolor=vthreecolor](81,12){\cirrad}
\psline[linecolor=vthreecolor](81,13)(81,14)
\pscircle*[linecolor=vthreecolor](81,13){\cirrad}
\pscircle*[linecolor=vthreecolor](81,14){\cirrad}
\pscircle*[linecolor=vzerocolor](82,17){\cirrad}
\psline[linecolor=vonecolor](82,2)(83,3)
\pscircle*[linecolor=vonecolor](82,2){\cirrad}
\psline[linecolor=vthreecolor](82,7)(82,8)
\psline[linecolor=vthreecolor](82,7)(85,9)
\pscircle*[linecolor=vthreecolor](82,7){\cirrad}
\psline[linecolor=vthreecolor](82,8)(82,9)
\pscircle*[linecolor=vthreecolor](82,8){\cirrad}
\pscircle*[linecolor=vthreecolor](82,9){\cirrad}
\psline[linecolor=vthreecolor](82,6)(85,6)
\pscircle*[linecolor=vthreecolor](82,6){\cirrad}
\psline[linecolor=vzerocolor](83,14)(83,15)
\psline[linecolor=vzerocolor](83,14)(84,15)
\psline[linecolor=vzerocolor](83,14)(86,15)
\pscircle*[linecolor=vzerocolor](83,14){\cirrad}
\psline[linecolor=vonecolor](83,1)(83,2)
\pscircle*[linecolor=vonecolor](83,1){\cirrad}
\psline[linecolor=vonecolor](83,2)(83,3)
\pscircle*[linecolor=vonecolor](83,2){\cirrad}
\pscircle*[linecolor=vonecolor](83,3){\cirrad}
\psline[linecolor=vthreecolor](83,10)(86,10)
\pscircle*[linecolor=vthreecolor](83,10){\cirrad}
\pscircle*[linecolor=vthreecolor](83,18){\cirrad}
\psline[linecolor=vthreecolor](83,15)(83,16)
\psline[linecolor=vthreecolor](83,15)(86,16)
\pscircle*[linecolor=vthreecolor](83,15){\cirrad}
\psline[linecolor=vthreecolor](83,16)(86,17)
\pscircle*[linecolor=vthreecolor](83,16){\cirrad}
\pscircle*[linecolor=vzerocolor](84,22){\cirrad}
\pscircle*[linecolor=vthreecolor](84,19){\cirrad}
\pscircle*[linecolor=vthreecolor](84,12){\cirrad}
\psline[linecolor=vthreecolor](84,5)(85,6)
\pscircle*[linecolor=vthreecolor](84,5){\cirrad}
\psline[linecolor=vthreecolor](84,15)(85,16)
\pscircle*[linecolor=vthreecolor](84,15){\cirrad}
\psline[linecolor=vzerocolor](85,20)(86,21)
\pscircle*[linecolor=vthreecolor](85,20){\cirrad}
\psline[linecolor=vthreecolor](85,8)(85,9)
\psline[linecolor=vthreecolor](85,8)(86,10)
\pscircle*[linecolor=vthreecolor](85,8){\cirrad}
\pscircle*[linecolor=vthreecolor](85,9){\cirrad}
\psline[linecolor=vthreecolor](85,19)(88,19)
\pscircle*[linecolor=vthreecolor](85,19){\cirrad}
\pscircle*[linecolor=vthreecolor](85,6){\cirrad}
\pscircle*[linecolor=vthreecolor](85,13){\cirrad}
\psline[linecolor=vthreecolor](85,16)(86,17)
\pscircle*[linecolor=vthreecolor](85,16){\cirrad}
\pscircle*[linecolor=vzerocolor](86,21){\cirrad}
\psline[linecolor=vthreecolor](86,11)(86,12)
\psline[linecolor=vthreecolor](86,11)(87,14)
\pscircle*[linecolor=vthreecolor](86,11){\cirrad}
\psline[linecolor=vthreecolor](86,12)(86,13)
\pscircle*[linecolor=vthreecolor](86,12){\cirrad}
\pscircle*[linecolor=vthreecolor](86,13){\cirrad}
\pscircle*[linecolor=vthreecolor](86,10){\cirrad}
\pscircle*[linecolor=vthreecolor](86,20){\cirrad}
\psline[linecolor=vthreecolor](86,15)(86,16)
\psline[linecolor=vthreecolor](86,15)(89,15)
\pscircle*[linecolor=vthreecolor](86,15){\cirrad}
\psline[linecolor=vthreecolor](86,16)(86,17)
\pscircle*[linecolor=vthreecolor](86,16){\cirrad}
\pscircle*[linecolor=vthreecolor](86,17){\cirrad}
\pscircle*[linecolor=vthreecolor](86,8){\cirrad}
\psline[linecolor=vonecolor](87,2)(87,3)
\pscircle*[linecolor=vonecolor](87,2){\cirrad}
\psline[linecolor=vonecolor](87,3)(87,4)
\pscircle*[linecolor=vonecolor](87,3){\cirrad}
\pscircle*[linecolor=vonecolor](87,4){\cirrad}
\pscircle*[linecolor=vthreecolor](87,14){\cirrad}
\pscircle*[linecolor=vthreecolor](87,6){\cirrad}
\psline[linecolor=vonecolor](88,2)(89,3)
\pscircle*[linecolor=vonecolor](88,2){\cirrad}
\psline[linecolor=vthreecolor](88,8)(89,9)
\pscircle*[linecolor=vthreecolor](88,8){\cirrad}
\psline[linecolor=vthreecolor](88,18)(88,19)
\pscircle*[linecolor=vthreecolor](88,18){\cirrad}
\pscircle*[linecolor=vthreecolor](88,19){\cirrad}
\pscircle*[linecolor=vthreecolor](88,10){\cirrad}
\pscircle*[linecolor=vthreecolor](88,12){\cirrad}
\pscircle*[linecolor=vonecolor](89,3){\cirrad}
\pscircle*[linecolor=vthreecolor](89,9){\cirrad}
\pscircle*[linecolor=vthreecolor](89,15){\cirrad}
\psline[linecolor=vthreecolor](89,7)(90,8)
\pscircle*[linecolor=vthreecolor](89,7){\cirrad}
\pscircle*[linecolor=vonecolor](90,2){\cirrad}
\psline[linecolor=vthreecolor](90,6)(90,7)
\pscircle*[linecolor=vthreecolor](90,6){\cirrad}
\psline[linecolor=vthreecolor](90,7)(90,8)
\pscircle*[linecolor=vthreecolor](90,7){\cirrad}
\pscircle*[linecolor=vthreecolor](90,8){\cirrad}
\end{pspicture}

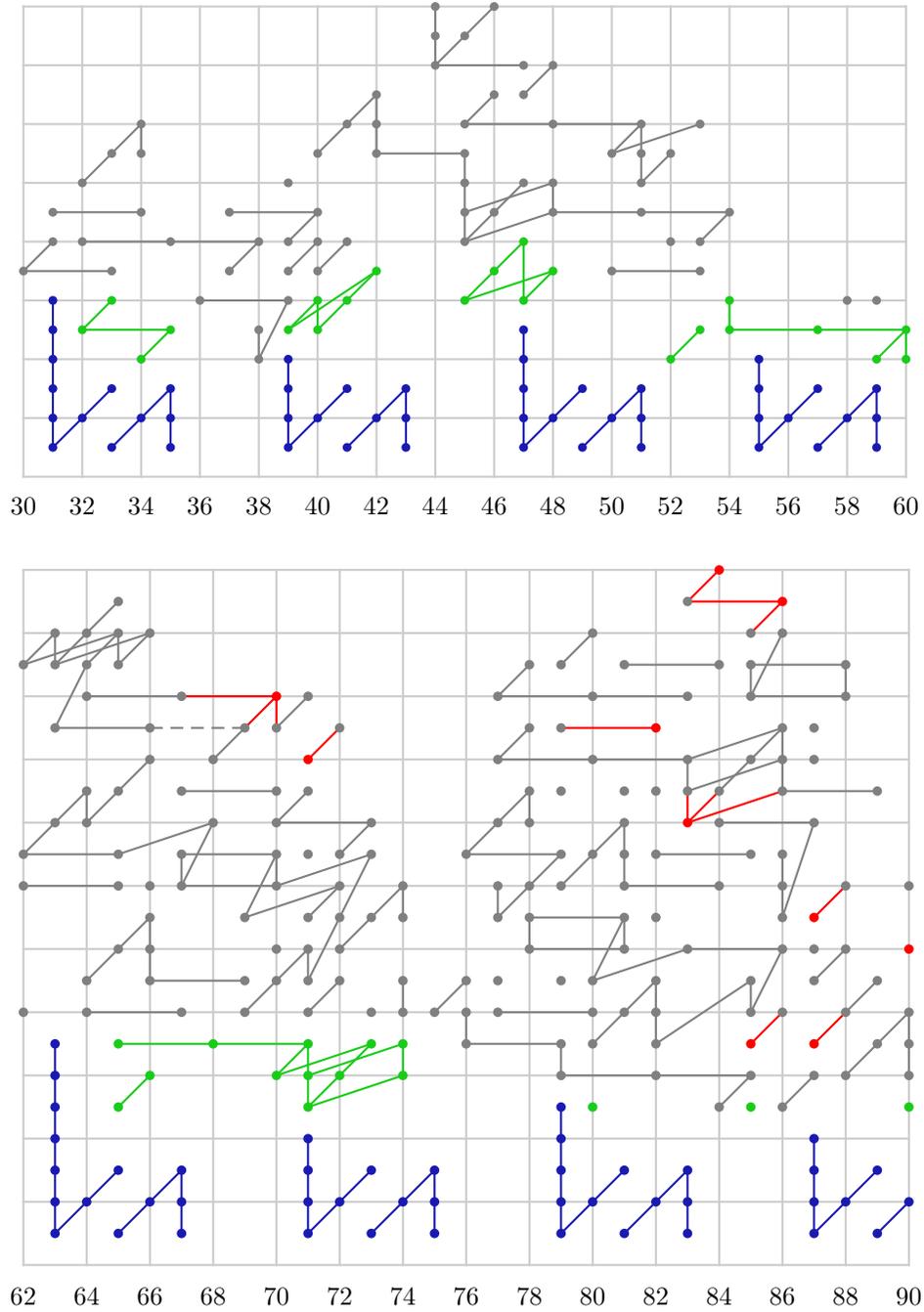
Figure \ref{fig:Hatcher} displays the $2$-primary stable homotopy groups 
in a graphical format.  Vertical chains of $n$ dots in column $k$
indicate a copy of $\Z/2^n$ in $\pi_k^{\wedge}$.  The non-vertical
lines indicate additional multiplicative structure that we will not
discuss here.  The blue dots represent the $v_1$-periodic subgroups.
The green dots are associated to the topological modular forms
spectrum $\mathit{tmf}$; the precise relationship is too complicated
to describe here.  Finally, the red dots indicate uncertainties.
Allen Hatcher originally promoted this type of graphical description 
of stable homotopy groups.

The orders of individual $2$-primary stable homotopy groups do not 
follow a clear pattern, with large increases and decreases seemingly
at random.  However, an empirically observed pattern emerges
if we consider the cumulative size of the groups, i.e., the product
of the orders of 
all $2$-primary stable homotopy groups from dimension $1$ to
dimension $k$.  

Our data strongly suggests that asymptotically, there is a linear relationship between $k^2$ 
and the logarithm of this product of orders.
In other words, the number of dots in Figure \ref{fig:Hatcher} in
stems $1$ through $k$ is linearly proportional to $k^2$.
Thus, in extending from dimension 60 to dimension 90, the overall 
size of the computation more than doubles.

\begin{conj}
Let $f(k)$ be the product of the orders of the $2$-primary
stable homotopy groups in dimensions $1$ through $k$.
Then $\log f(k) = O(k^2)$.
\end{conj}

One interpretation of this conjecture is that the expected value of
the logarithm of the order of the 2-primary component of $\pi_k$
grows linearly in $k$.

\section{$v_1$-periodic stable homotopy groups}
\label{sctn:v1}

Adams provided the first
infinite families of elements in the stable homotopy groups \cite{Adams66}.
Within the stable homotopy groups, there is a regular repeating pattern
of subgroups.
These subgroups are
known as the ``$v_1$-periodic stable homotopy
groups", and they are closely related to the ``image of $J$".

In order to describe the $2$-primary component of the
$v_1$-periodic stable homotopy groups precisely,
we will need the following elementary number-theoretic definition.
For any integer $k$, write $k$ as a product $2^a j$, where
$j$ is odd.  Define $\nu(k)$ to be $2^a$.

\begin{thm} \cite{Adams66}
\label{thm:v1}
Table \ref{tab:v1} gives the $v_1$-periodic stable homotopy groups
inside of $\pi_k^\wedge$ for all $k \geq 2$.
\end{thm}

\begin{longtable}{l|l}
\caption{$2$-primary $v_1$-periodic stable homotopy groups
\label{tab:v1} 
} \\
\toprule
$k$ & $v_1$-periodic   \\
& subgroup of $\pi_k^\wedge$ \\
\midrule \endfirsthead
\caption[]{$2$-primary $v_1$-periodic stable homotopy groups} \\
\toprule
$k$ & $v_1$-periodic  \\
 & subgroup of $\pi_k^\wedge$ \\
\midrule \endhead
\bottomrule \endfoot
$0 \mod 8$ & $2$ \\
$1 \mod 8$ & $2^2$ \\
$2 \mod 8$ & $2$ \\
$3 \mod 8$ & $8$ \\
$7 \mod 8$ & $2 \nu(k+1)$ \\
\end{longtable}

\begin{remark}
For an odd prime $p$, 
one can also easily describe the $p$-primary $v_1$-periodic
stable homotopy groups.  The groups are zero unless 
$k \equiv -1 \mod 2(p-1)$.  
If $k + 1 = 2(p-1) a$, then let $\nu_p(a)$ be the largest factor of $a$
that is a power of $p$.  The $p$-primary $v_1$-periodic
stable homotopy group is a cyclic group of order $p \cdot \nu_p(a)$.
\end{remark}

The $v_1$-periodic subgroups are direct summands of
the stable homotopy groups.  Their complementary
summands are known as the ``$v_1$-torsion subgroups".
The language involving periodicity and torsion derives from
the theory of nilpotence and periodicity due
to Devinatz, Hopkins, Ravenel, Smith, and others
\cite{DHS88} \cite{HS98} \cite{Ravenel92},
which goes beyond the scope of this article.

\section{The Adams spectral sequence}
\label{sctn:Adams}

The most powerful tool for computing $\pi_k^\wedge$ is
the Adams spectral sequence \cite{Adams58}.
Information about
ordinary cohomology, together with its higher structure in the
form of cohomology operations, yields information about the
stable homotopy groups.
See \cite{IWX20b} for a graphical representation of the Adams
spectral sequence.

The Steenrod algebra $A$ is the ring of stable operations on cohomology
with $\F_2$ coefficients \cite{Steenrod62}.  It is generated by the 
Steenrod squaring operations $\Sq^i$ of degree $i$.  
It is convenient to adopt the convention that $\Sq^0 = 1$.
These
operations satisfy the Adem relations
\[
\Sq^a \Sq^ b = \sum_{c = 0}^{a+b} \binom{b-c-1}{a-2c} \Sq^{a+b-c} \Sq^c
\]
whenever $a < 2b$.
The Steenrod algebra also has a ``coproduct" that takes
$\Sq^i$ to $\sum_{c = 0}^i \Sq^{i-c} \otimes \Sq^c$.
The product and coproduct structure on the Steenrod algebra
together form the structure of a Hopf algebra.

Unfortunately, the Adem relations are a bit hard to grasp.  
Following ideas of Milnor \cite{Milnor58}, it turns out to be much easier
to work with the $\F_2$-dual.  In other words, we consider
$A_* = \Hom_{\F_2} (A, \F_2)$.  Then the product and coproduct on
$A$ become a coproduct and a product on $A_*$ respectively,
and $A_*$ is again a Hopf algebra.

The Hopf algebra $A_*$ is much easier to describe.
It is isomorphic to the polynomial ring
$\F_2 [\zeta_1, \zeta_2, \ldots ]$, where the coproduct is given
by the formula
\[
\zeta_i \mapsto \sum_{k = 0}^i \zeta_{i-k}^{2^k} \otimes \zeta_{k}.
\]
By convention, we let $\zeta_0$ equal $1$.
The duality between the structure on $A_*$ and on $A$
is not obvious.

The next step in the Adams spectral sequence program is to consider
the derived groups $\Ext_A(\F_2, \F_2)$ of 
$\Hom_A(\F_2, \F_2)$, in the sense of homological algebra.
Roughly speaking, these $\Ext$ groups capture higher information
about the structure of $A$, including generators, relations, 
relations among relations, etc.

The groups $\Ext_A(\F_2, \F_2)$ serve as the input to the 
Adams spectral sequence.  Up to dimension $13$, these algebraic
$\Ext$ groups give complete information about $\pi_k^\wedge$.
However, in higher dimensions, further complications occur.
Specifically, one must compute Adams differentials.  These differentials
measure the deviation between algebra and homotopy.  In practice,
the computation of these differentials is the limiting factor
in our knowledge of stable homotopy groups.

In higher dimensions, the algebraic groups $\Ext_A(\F_2,\F_2)$
themselves become difficult to compute directly.
The May spectral sequence  
is the best way to compute these groups by hand \cite{May64}.
With the aid of this spectral sequence, May extended the
computations of stable homotopy groups to approximately
dimension 30.
Tangora eventually carried out the algebraic computation of
$\Ext_A(\F_2, \F_2)$ via the May spectral sequence to
dimension 70 \cite{Tangora70a}.

In the modern era, the most efficient way to compute $\Ext$ groups
is by machine. 
Bruner \cite{Bruner89} \cite{Bruner93} \cite{Bruner97}, Nassau  \cite{Nassau}, and the second author 
\cite{Wang19} have constructed various
efficient algorithms that provide a wealth of algebraic data,
far surpassing our ability to interpret.  
The most extensive
computations extend beyond dimension 200.
For practical purposes, we can take this $\Ext$ data as given.
Computer assisted techniques are likely 
to continue to grow in importance in the
computation of stable homotopy groups.

Beyond dimension 30, the analysis of Adams differentials becomes
more difficult \cite{BMT70} \cite{Bruner84} \cite{Isaksen14}
\cite{MT67}.  
The stable homotopy groups possess higher structure in the form
of Massey products and Toda brackets, and this higher structure
leads to additional information about differentials.
However, these arguments are notoriously tricky, 
and the published literature contains
more than one example of an incorrect computation.

The practical limit of this style of argument occurs at 
dimension 61.  See \cite{Isaksen14} for a thorough accounting
of the Adams spectral sequence through dimension 59.
The article \cite{WX17} employs strenuous efforts to obtain
just two more stable homotopy groups in dimensions 60 and 61.
Beyond dimension 61, these methods are simply no longer practical.

\section{The motivic Adams spectral sequence}

Morel and Voevodsky \cite{Morel99a} \cite{MV99}  
developed motivic homotopy theory
in the mid 1990's as a means of importing homotopical techniques
into algebraic geometry.  This program found great success in 
Voevodsky's
resolutions of the Milnor Conjecture \cite{Voevodsky03b} and the
Bloch-Kato Conjecture \cite{Voevodsky11}.
For our purposes, we may simplify the theory somewhat by
considering only cellular objects and by taking appropriate completions
at a prime $p$.

Motivic homotopy theory is bigraded, so all invariants,
including cohomology and stable homotopy groups, are bigraded.
There is a bigraded family of spheres $S^{p,q}$ that serve
as the basic building blocks of motivic homotopy theory.
While motivic homotopy theory can be studied over any base field,
we will focus only on the case when the base field is $\C$.
The additional structure contained in
$\C$-motivic homotopy theory provides a new tool for computing
classical stable homotopy groups
\cite{IWX}.

The use of $\C$-motivic homotopy theory suggests that these 
stable homotopy groups computations
are logically dependent on deep and difficult algebro-geometric
results of Voevodsky on the motivic cohomology of a point 
\cite{Voevodsky03b}
and the structure of the motivic Steenrod algebra 
\cite{Voevodsky03a} \cite{Voevodsky10}.
However, 
there are now two entirely topological models
for the part of $\C$-motivic homotopy theory that is relevant
to stable homotopy group computations \cite{GIKR18} \cite{Pstragowski18}.
The fundamental inputs are explicit descriptions of
the cohomology of a point and of the Steenrod algebra,
and these inputs can be
derived from first principles in these topological models, using
nothing more than well-known standard classical computations.
Therefore, our new stable homotopy group computations are not
actually
logically dependent on anything algebro-geometric.

The $\C$-motivic cohomology of a point takes the form
$\F_2[\tau]$, where $\tau$ has degree $(0,1)$.
The dual $\C$-motivic Steenrod algebra $A^\C_*$ takes the form
\[
\frac{\F_2[\tau][\tau_0, \tau_1, \ldots, \xi_1, \xi_2, \ldots]}
{\tau_i^2 = \tau \xi_{i+1}},
\]
where the coproduct is given by the formulas
\[
\tau_i \mapsto \tau_i \otimes 1 +
\sum_{k = 0}^i \xi_{i-k}^{2^k} \otimes \tau_{k}
\]
\[
\xi_i \mapsto 
\sum_{k = 0}^i \xi_{i-k}^{2^k} \otimes \xi_{k}.
\]
By convention, we let $\xi_0$ equal $1$.

Comparison to the classical dual Steenrod algebra illuminates
the subtleties of the $\C$-motivic dual Steenrod algebra.
After inverting $\tau$, the element $\xi_{i+1}$ becomes decomposable,
so $A^{\C}_*[\tau^{-1}]$ is a polynomial algebra
over $\F_2[\tau^{\pm 1}]$ on generators $\tau_i$.
This recovers the classical dual Steenrod algebra, where
$\tau_i$ and $\xi_{i+1}$ correspond to $\zeta_i$ and $\zeta_i^2$
respectively.

On the other hand, after setting $\tau$ equal to zero,
the result is an exterior algebra on generators $\tau_i$
tensored with a polynomial algebra on generators $\xi_{i+1}$.
This structure is analogous to Milnor's description of the
 classical dual Steenrod algebra
at odd primes \cite{Milnor58}.

As in the classical case,
$\Ext_{A^\C}(\F_2[\tau], \F_2[\tau])$ can be computed
by machine in a large range.  Then  
$\C$-motivic Adams differentials can be determined 
by the standard methods.
As in the classical case, dimension 61 seems to be the 
practical limit of this approach.

\section{Algebraicity of the cofiber of $\tau$}

There is a map $S^{0,-1} \map S^{0,0}$ in the 
$\C$-motivic stable homotopy category that induces multiplication
by $\tau$ in $\C$-motivic cohomology.  Therefore, we use the same
notation $\tau$ for this map between spheres.

Let $S/\tau$ be the mapping cone (or cofiber) of $\tau$.
This object is a stable 2-cell complex that a priori
has no special structure.  Surprisingly, the homotopically 
defined $S/\tau$ has a remarkably algebraic structure.

Recall that $BP$ is the Brown-Peterson generalized cohomology
theory \cite{BP66}.  This cohomology theory
has been of remarkable use in the computation of stable
homotopy groups \cite{Ravenel86}.

The ring
of coefficients for this theory is
$BP_* = \Z_{(2)}[v_1, v_2, \ldots ]$, and 
$BP_* BP$ is the dual of the ring of stable operations.
The full structure of the object $BP_* BP$ can be completely
described, although we do not give the details here.
Then 
$\Ext_{BP_* BP}(BP_*, BP_*)$ is the $E_2$-page of the
Adams-Novikov spectral sequence, which is another tool
for computing stable homotopy groups that is complementary
to the Adams spectral sequence.
These $\Ext$ groups themselves are quite complicated, but
they can be computed in a range by machine.
Alternatively, they can be computed by the algebraic Novikov
spectral sequence
\cite{Miller75} \cite{Novikov67}.

\begin{thm}
\label{thm:S/t}
\cite{GWX18}
The $\C$-motivic Adams spectral sequence that computes the 
motivic stable homotopy groups of $S/\tau$ is isomorphic to
the algebraic Novikov spectral sequence.
\end{thm}

In fact, Theorem \ref{thm:S/t} is a computational corollary of
other more structural
results.  In particular, Gheorghe demonstrated that
the $\C$-motivic spectrum $S/\tau$ is an $E_\infty$-ring 
object in an essentially unique way \cite{Gheorghe17}, and 
the homotopy category of cellular $S/\tau$-modules is equivalent
to a derived category of $BP_* BP$-comodules
\cite{GWX18}.

Deformation theory provides a unifying perspective on this circle of ideas.
The key insight is that $\C$-motivic 
cellular stable homotopy theory is a deformation of classical stable
homotopy theory \cite{GWX18}, after completing at each prime $p$.  
From this perspective, the ``generic fiber" of
$\C$-motivic stable homotopy theory is classical stable homotopy theory,
and the ``special fiber"
has an entirely algebraic description.  The special fiber is
 the category of $BP_*BP$-comodules, or equivalently,
the category of quasicoherent sheaves on the moduli stack of 
1-dimensional formal groups.  

Theorem \ref{thm:S/t} is
particularly valuable for computation.  It means that
the Adams spectral sequence for $S/\tau$
can be computed in an entirely algebraic manner, i.e.,
can be computed by machine in a large range.
This observation leads to the following innovative
program for computing classical stable homotopy groups.

\begin{enumerate}
\item
\label{item:machine-Adams}
Compute the
$\C$-motivic Adams $E_2$-page by machine in a large range.
\item
\label{item:machine-algNov}
Compute
the algebraic Novikov spectral sequence by machine in a large
range, including
all differentials and multiplicative structure.
\item
Use Theorem \ref{thm:S/t}
to deduce the structure of the
motivic Adams spectral sequence for $S/\tau$.
\item
Use the cofiber sequence
\[
S^{0,-1} \stackrel{\tau}{\to} S^{0,0} \to S/\tau \to 
S^{1,-1}
\]
and naturality of Adams spectral sequences 
to pull back and push forward Adams differentials for $S/\tau$
to Adams differentials for the motivic sphere.
\item
\label{item:ad-hoc}
Apply a variety of ad hoc arguments to deduce additional Adams differentials
for the motivic sphere.
\item
Use a long exact sequence in homotopy groups to deduce hidden $\tau$
extensions in the motivic Adams spectral sequence for the sphere.
\item
Invert $\tau$ to obtain the classical Adams spectral sequence and the 
classical stable homotopy groups.
\end{enumerate}

The machine-generated data that we use in steps (\ref{item:machine-Adams})
and (\ref{item:machine-algNov}) are available at \cite{Wang19}.

As the dimension increases, the ad hoc arguments of step
(\ref{item:ad-hoc}) become more and more complicated.
Eventually, this approach will break down when the ad hoc
arguments become too complicated to resolve.
It is not yet clear when that will occur.

\section{Remaining uncertainties}
\label{sctn:unknown}

Up to dimension 90, there are only four Adams differentials
whose values have not been completely determined.
Each of these unknown differentials leads to uncertainties
in the stable homotopy groups in two adjacent dimensions.
The first such differential affects the orders of
$\pi_{69}^\wedge$ and $\pi_{70}^\wedge$.
In most cases, these uncertainties mean that the order of $\pi_k$ is 
known only up to a factor of $2$.  In a few cases, the additive group
structure is also undetermined.
These uncertainties are not independent, as described in the following
alternatives.

In dimensions 70 and 71, one of the following two possibilities
occurs:
\begin{enumerate}
\item
the $v_1$-torsion in $\pi_{70}^\wedge$ has order $512$;
and the $v_1$-torsion in $\pi_{71}^\wedge$ has order $2048$.
\item
the $v_1$-torsion in $\pi_{70}^\wedge$ has order $256$;
and the $v_1$-torsion in $\pi_{71}^\wedge$ has order $1024$.
\end{enumerate}
In dimensions 82 and 83, one of the following two possibilities occurs:
\begin{enumerate}
\item
the $v_1$-torsion in $\pi_{82}^\wedge$ has order $256$, and
the $v_1$-torsion in $\pi_{83}^\wedge$ has order $64$.
\item
the $v_1$-torsion in $\pi_{82}^\wedge$ has order $128$, and
the $v_1$-torsion in $\pi_{83}^\wedge$ has order $32$.
\end{enumerate}
In dimensions 84 and 85, one of the following two possibilities
occurs:
\begin{enumerate}
\item
the $v_1$-torsion in $\pi_{84}^\wedge$ has order $64$, and
the $v_1$-torsion in $\pi_{85}^\wedge$ has order $1024$.
\item
the $v_1$-torsion in $\pi_{84}^\wedge$ has order $32$, and
the $v_1$-torsion in $\pi_{85}^\wedge$ has order $512$.
\end{enumerate}
In dimensions $86$ and $87$, one of the following two possibilities
occurs:
\begin{enumerate}
\item
the $v_1$-torsion in $\pi_{86}^\wedge$ has order $2048$, and
the $v_1$-torsion in $\pi_{87}^\wedge$ has order $256$.
\item
the $v_1$-torsion in $\pi_{86}^\wedge$ has order $1024$, and
the $v_1$-torsion in $\pi_{87}^\wedge$ has order $128$.
\end{enumerate}

\bibliographystyle{amsalpha}

\end{document}